\documentclass[leqno]{amsart}
\usepackage{amssymb,amsthm,amsmath,color,dsfont}
\usepackage{varioref} 
\usepackage[square,numbers]{natbib}
\usepackage{hyperref}  
\usepackage{aliascnt}
\usepackage{graphicx}

\numberwithin{equation}{section}

    \labelformat{equation}{\textup{(#1)}}

    \newtheorem{thm}{Theorem}[section]

    \newaliascnt{fact}{thm}
    
    \aliascntresetthe{fact}

    \newaliascnt{conj}{thm}
    \newtheorem{conj}[conj]{Conjecture}
    \aliascntresetthe{conj}

    \newaliascnt{lem}{thm}
    \newtheorem{lem}[lem]{Lemma}
    \aliascntresetthe{lem}

    \newaliascnt{cor}{thm}
    
    \aliascntresetthe{cor}

    \newaliascnt{prop}{thm}
    
    \aliascntresetthe{prop}

    \theoremstyle{definition}
    \newaliascnt{defn}{thm}
    \newtheorem{defn}[defn]{Definition}
    \aliascntresetthe{defn}

    \theoremstyle{remark}
    \newaliascnt{rem}{thm}
    \newtheorem{rem}[rem]{Remark}
    \aliascntresetthe{rem}

\renewcommand{\P}{\mathbb P}
\newcommand{\E}{\mathbb E}
\newcommand{\si}{\sigma}
\usepackage{dsfont}

\def\Z{\mathbb Z}
\begin{document}

\title{Central Binomial Tail Bounds}

\author{Matus Telgarsky}
\address{Department of Computer Science and Engineering, University of California, San Diego, 92093}
\email{mtelgars@cs.ucsd.edu}





\begin{abstract}
An alternate form for the binomial tail is presented, which 
leads to a variety of bounds for
the central tail.  A few can be weakened into the corresponding 
Chernoff and Slud bounds, which not only demonstrates the quality of the presented bounds, but
also provides alternate proofs for the classical bounds.
\end{abstract}

\maketitle
\section{Introduction}
Let $B(p,n)$ denote a binomial random variable comprising $n$ flips of
a bias-$p$ coin, and set $\si = \sqrt{p(1-p)}$.  The classical form of the
central tail, obtained by summing over the possible outcomes, is
\begin{equation}
\label{eq:classical_tail}
\P[B(p,n) \geq n/2] = \sum_{h=\lceil n/2\rceil}^n \binom{n}{h}p^h(1-p)^{n-h}.
\end{equation}
If instead $B$ is considered from the perspective of random walks on the
integer line, another representation is possible by tracking, as $n$ increases,
the motion of mass from one side of the origin to the other.  This intuition
is formalized in \autoref{sec:basics}, and results in the following statement.
\begin{thm}
\label{thm:ctail}
When $n$ is odd and $p < 1/2$,
\begin{align}
\P[B(p,n) \geq n/2] 
&= p - (1/2-p) \sum_{j=1}^{(n-1)/2} \binom{2j}{j}\si^{2j} \label{eq:ctail_finite}\\
&= (1/2-p) \sum_{j\geq (n+1)/2} \binom{2j}{j} \si^{2j}. \label{eq:ctail_infinite}
\end{align}
\end{thm}
To demonstrate the value of this new characterization of the central tail, it
is used to derive bounds.  In particular, \autoref{sec:elembounds}
(``Closed-form Bounds'') 
approximates the summands of \ref{eq:ctail_finite}
and \ref{eq:ctail_infinite} in various ways to yield summations with
closed-form expressions.  On the other hand, \autoref{sec:contbounds}
(``Bounding with the Standard Normal'') 
replaces the summation of \ref{eq:ctail_infinite} with an integral,
yielding a bound incorporating the distribution function of the standard normal.

\begin{figure}[t]
\includegraphics[width=\textwidth]{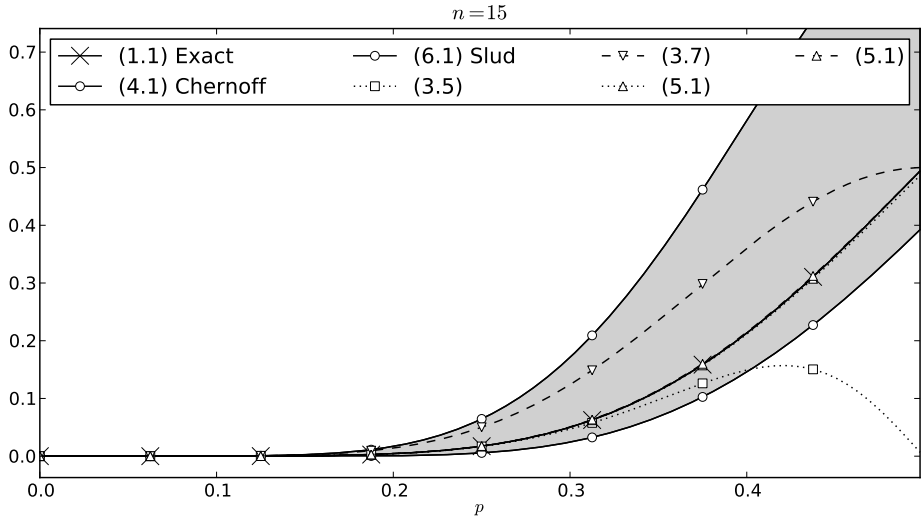}
\caption{
Bounds are identified by equation number, lower and upper bounds derived in this paper being
respectively dotted and dashed.  Shading identifies
the region between the Chernoff and Slud bounds.
The exact (central) binomial tail is
marked with x's, however the bounds in \ref{eq:cont_bounds} somewhat occlude
it.}
\label{fig:smallplot}
\end{figure}

A few of these bounds appear in \autoref{fig:smallplot}.  The upper and lower bound 
pair of \ref{eq:cont_bounds} are the tightest, and their forms are sufficiently similar
to allow the gap to be analytically quantified.  This comes at the cost of interpretability:
they are rather complicated.  Contrastingly, the upper bound of
\ref{eq:ctail_chernlike} is simple: when $n$ is odd, it is just
$(2\si)^{n+1}/2$.  Remarkably, the Chernoff bound for $m$ even is $(2\si)^m$.
(Whenever the number of trials is odd, $n$ is used; when it is even, $m$ is
used.)  The relationship between the bounds of \autoref{sec:elembounds} and the Chernoff
bound is explored in \autoref{sec:chern}. The gap between these bounds,
seemingly large in \autoref{fig:smallplot}, is quantified, which furthermore
provides an alternate proof of the Chernoff bound.

Among the lower bounds, only that of \ref{eq:cont_bounds} consistently outperforms Slud's bound.
Fortunately, both depend on the standard normal, and thus the comparison can be made precise:
\autoref{sec:slud} discusses proving Slud's bound by weakening the bounds of \ref{eq:cont_bounds}.
 The proof is notable because it extends the sufficient conditions of
the classical statement of Slud's bound.  Unfortunately, the details of this proof are tedious, and relegated to
\autoref{sec:pf_slud}.
The task of producing a good, elementary
lower bound proved challenging; the lower bound of
\ref{eq:funky_exp_lower_bound}, which appears in \autoref{fig:smallplot}, is
only tight for $p$ away from $1/2$.  A comparison of all bounds may be found in \autoref{fig:plots}
of \autoref{sec:allbounds}.

To close, \autoref{sec:gtail} generalizes \autoref{thm:ctail} to arbitrary tails,
however no bounds are derived.

\section{Central Binomial Tails via Random Walks}
\label{sec:basics}
Consider random walks on the integer line originating at 0, and at each step incrementing 
their position with probability $p$, or decrementing it with probability $1-p$.  When $n$ is
odd, $\P[B(p,n)\geq n/2]$ can be interpreted to mean the probability mass of walks terminating
with positive coordinate after $n$ steps.  To prove \autoref{thm:ctail}, the
first step is to
quantify the effect which two trials have on this probability mass.
\begin{lem}
\label{lem:ctail_step}
When $n$ is odd,
\[
\P[B(p,n+2)\geq (n+2)/2] - \P[B(p,n) \geq n/2] = (p-1/2)\binom{n+1}{(n+1)/2}\si^{n+1}.
\]
\end{lem}
\begin{proof}
For mass to change sign in two steps, it must originate in a path ending at
a coordinate adjacent to the origin, and move in the direction of the origin
twice.  Symbolically,
\begin{align*}
&\underbrace{p^2\binom{n}{(n-1)/2}p^{(n-1)/2}(1-p)^{(n+1)/2}}_{\textup{two increasing steps}}
-\underbrace{(1-p)^2\binom{n}{(n+1)/2}p^{(n+1)/2}(1-p)^{(n-1)/2}}_{\textup{two decreasing steps}}\\
&\quad\quad= (2p-1)\binom{n}{(n+1)/2}\si^{n+1}.
\end{align*}
To finish, note that 
$\binom{n}{(n+1)/2} = \frac {(n+1)/2}{n+1}\binom{n+1}{(n+1)/2} = \frac 1 2 \binom{n+1}{(n+1)/2}$.
\end{proof}
Although the above proof depends on a random walk interpretation, it also goes through
purely algebraically using \ref{eq:classical_tail}. 

Accumulating the contribution of such steps up to $n$, 
it is possible to rewrite the binomial tail.
\begin{lem}
\label{lem:ctail_finite}
When $n$ is odd,
\[
\P[B(p,n)\geq n/2] = p - (1/2-p) \sum_{j=1}^{(n-1)/2} \binom{2j}{j}\si^{2j}.
\]
\end{lem}
\begin{proof}Invoking \autoref{lem:ctail_step},
\begin{align*}
\P[B(p,n)\geq n/2] & \\
& \hspace{-0.75in}=\P[B(p,1) \geq 1/2]
+ \sum_{\substack{j=1\\ j\textup{ odd}}}^{n-2} \left(\P[B(p,j+2)\geq (j+2)/2] - \P[B(p,j)\geq j/2]\right)\\
& \hspace{-0.75in}=p
+(p-1/2) \sum_{\substack{j=1\\ j\textup{ odd}}}^{n-2} \binom{j+1}{(j+1)/2}\si^{j+1},
\end{align*}
and substituting $2j-1$ for $j$ yields the lemma.
\end{proof}

When $p<1/2$, as $n\to\infty$, the central tail probability must approach $0$ (cf. for instance \ref{eq:ctail_chernlike}).
As such, it should also be possible to compute the tail in a fashion complementary to
\autoref{lem:ctail_finite}, instead accumulating the
contribution of all remaining steps.

\begin{lem}
\label{lem:ctail_infinite}
When $n$ is odd and $p < 1/2$,
\[
\P[B(p,n)\geq n/2] = (1/2-p) \sum_{j\geq (n+1)/2} \binom{2j}{j} \si^{2j}.
\]
\end{lem}
\begin{proof}
When $p=0$, the result is immediate, thus take $p\in(0,1/2)$.  Combining the Taylor expansion $(1-4\si^2)^{-1/2} = \sum_{j\geq0} \binom{2j}{j}\si^{2j}$ with
\autoref{lem:ctail_finite},
\begin{align*}
p - (\frac 1 2-p) \sum_{j=1}^{(n-1)/2} \binom{2j}{j}\si^{2j}
&= p - (\frac 1 2-p)  \left(\frac 1 {\sqrt{1-4\si^2}} - \sum_{j\geq(n+1)/2} \binom{2j}{j}\si^{2j} - 1\right)\\
&= \frac 1 2 - \frac {1/2-p} {\sqrt{1-4\si^2}} + (1/2-p)\sum_{j\geq(n+1)/2} \binom{2j}{j}\si^{2j},
\end{align*}
the result following since $\sqrt{1-4\si^2} = |1-2p|$ and $p<1/2$.
\end{proof}

\begin{proof}[Proof of \autoref{thm:ctail}.] \ref{eq:ctail_finite} is handled by 
\autoref{lem:ctail_finite},
whereas \autoref{lem:ctail_infinite} takes
care of \ref{eq:ctail_infinite}.
\end{proof}

\begin{rem}
\label{rem:ctail_constraints}
Going forward, the two constraints that 
$n$ is odd
and $p < 1/2$ will frequently appear.
First note that $p<1/2$ can be assuaged with
\begin{align*}
\P[B(p,n)\geq n/2+k] 
&= \P[B(1-p,n) \leq n/2-k].
\end{align*}
(And when $p=1/2$, by \autoref{lem:ctail_finite}, $\P[B(1/2,n)\geq n/2] = 1/2$.)
Additionally, the tail may be flipped with
\[
\P[B(p,n)\geq n/2+k] 
= 1 - \P[B(p,n) < n/2+k].
\]
Lastly, using the same random walk reasoning as in the proof of
\autoref{lem:ctail_step}, central tail bounds on $B(p,m)$ where $m$ is even can be reduced to bounds on
$B(p,m-1)$ via
\begin{align*}
&\P[B(p,m)\geq m/2] - \P[B(p,m-1)\geq (m-1)/2] \\
&\quad\quad= p\binom{m-1}{(m-2)/2}p^{(m-2)/2}(1-p)^{m/2} 
= \frac 1 2\binom{m}{m/2}\si^{m}.
\end{align*}
\end{rem}

\section{Closed-form Bounds}
\label{sec:elembounds}
To produce bounds from \autoref{thm:ctail}, the first task is to eliminate the binomial coefficient,
which the following steps achieve by way of Stirling's approximation.
\begin{defn}
Define 
\begin{align*}
l(n) &= - \frac {9n+1}{3n(12n+1)},\\
u(n) &= - \frac {18n-1}{12n(6n+1)}.
\end{align*}
\end{defn}
\begin{rem}
Note that, for $n\geq 1$, both are strictly increasing, and 
\[
- \frac 1 {5n}
>
u(n)
>
- \frac 1 {4n}
>
l(n)
>
- \frac 1 {3n}
.
\]
\end{rem}
Using the bounded form of Stirling's 
approximation (as in (9.15) from chapter~2 of~\citet{feller})
\begin{equation*}
\sqrt{2\pi n} \left(\frac n e\right)^n e^{1/(12n+1)}
< n! < 
\sqrt{2\pi n} \left(\frac n e\right)^n e^{1/(12n)},
\end{equation*}
the central binomial coefficient can be bounded with
\begin{equation}
\label{eq:central_stirling}
\frac {4^je^{l(2j)}}{\sqrt{\pi j}} < \binom {2j}{j} < \frac {4^j e^{u(2j)}}{\sqrt {\pi j}}.
\end{equation}
Notice that combining \ref{eq:central_stirling} and (for instance) \autoref{lem:ctail_infinite}
yields the somewhat hopeful relation
\begin{equation}
\label{eq:ctail_stirling}
(1/2-p)\hspace{-0.10in}
\sum_{j\geq (n+1)/2} \hspace{-0.10in}
\frac {(2\si)^{2j}e^{l(2j)}} {\sqrt {\pi j}}
\leq
\P[B(p,n)\geq n/2]
\leq
(1/2-p)\hspace{-0.10in}
\sum_{j\geq (n+1)/2} \hspace{-0.10in}
\frac {(2\si)^{2j}e^{u(2j)}} {\sqrt {\pi j}}.
\end{equation}
The remainder of this section starts from \ref{eq:ctail_stirling} (or from
the analogous formula using the finite summation of \autoref{lem:ctail_finite}),
and manipulates the summation into one possessing a closed-form expression.
The primary difficulty in \ref{eq:ctail_stirling} is the term $j^{-1/2}$,
and the derivation of each bound can be characterized by its approach to this term.
A sense of the relative performance of the bounds can be gleaned from \autoref{fig:plots}
on page~\pageref{fig:plots}.

The first
bounds relax $j^{-1/2}$ trivially; that is, upper bounding it with 1, and lower bounding it with
$j^{-1}$.

\begin{thm}
When $n$ is odd
and $p\in(0,1/2)$,
\begin{align}
\P[B(p,n)\geq n/2] &\geq 
\frac {(1-2p)e^{l(n+1)}(2\si)^{n-1}}{\sqrt{2\pi(n+1)}}
\left(-\ln(1-4\si^2)
\right),\label{eq:elem_lb}\\
\P[B(p,n)\geq n/2] &\leq 
\frac {(2\si)^{n+1}}{(1-2p)\sqrt{2\pi(n+1)}}.
\label{eq:elem_ub}
\end{align}
\end{thm}
(Note that the lower and upper bounds may be related using $\ln(x) \leq x-1$.)
Both bounds become poor as $p\to1/2$; in fact, the upper bound
grows unboundedly, and the lower bound goes to zero.  The upper bound, however, is sufficiently
tight for $p<1/4$ and odd $n$ to prove the Chernoff bound in \autoref{sec:chern}.
\begin{proof}
From \ref{eq:ctail_stirling},
\begin{align*}
\P[B(p,n)\geq n/2]
&\geq(1/2-p) \sum_{j\geq (n+1)/2} \frac {(2\si)^{2j}e^{l(2j)}}{\sqrt{\pi j}}\\
&\geq
\frac {(1/2-p)e^{l(n+1)}}{\sqrt\pi} \sum_{j\geq 1} \frac {(2\si)^{2j+(n-1)}}{j\sqrt{(n+1)/2}},
\end{align*}
{\sloppy
and the 
bound follows using the Taylor expansion 
$-\ln(1-2\si^2) = \sum_{j\geq 1} (2\si)^{2j}/j$.}  Similarly for the upper bound,
\begin{align*}
\P[B(p,n)\geq n/2]
&\leq(1/2-p) \sum_{j\geq (n+1)/2} \frac {(2\si)^{2j}e^{u(2j)}}{\sqrt{\pi j}}\\
&\leq
\frac {1-2p}{\sqrt {2\pi(n+1)}} \sum_{j\geq (n+1)/2} (2\si)^{2j} \\
&=
\frac {1-2p}{\sqrt {2\pi(n+1)}}
 \left(\frac {(2\si)^{n+1}}{1-4\si^2}\right).
\end{align*}
To finish, use $1-4\si^2 = (1-2p)^2$.
\end{proof}

Another approach to the term $j^{-1/2}$ is to lower bound with an exponential.
\begin{thm}
\label{thm:funky_exp_lower_bound}
When $p<1/2$ and $n$ is odd,
\begin{equation}
\label{eq:funky_exp_lower_bound}
\P[B(p,n)\geq n/2] \geq 
\frac
{(1-2p)e^{l(n+1)}(2\si)^{n+1}}
{\left(1-4\si^2\sqrt{(n+1)/(n+3)}\right)\sqrt{2\pi(n+1)}}.
\end{equation}
\end{thm}
This lower bound also approaches zero as $p\to 1/2$, but is otherwise the tightest
elementary lower bound in this paper.  It can be seen to dominate the lower bound
of \ref{eq:elem_lb} by taking $n$ large and using a tangent approximation to $\ln$.
\begin{proof}
Fitting an exponential to $((n+1)/2)^{-1/2}$ and $((n+3)/2)^{-1/2}$ yields
\[
\frac 1 {\sqrt{(n+1)/2}} \left(\frac {\sqrt{(n+1)/2}}{\sqrt{(n+3)/2}}\right)^{j-{(n+1)/2}}
=
\sqrt{\frac 2 {n+1}} \left(\sqrt{\frac {n+1}{n+3}}\right)^{j-{(n+1)/2}}.
\]
Thus, again using \ref{eq:central_stirling} and \autoref{lem:ctail_infinite},
\begin{align*}
\P[B(p,n)\geq n/2]
&\geq
\frac {(1/2-p)e^{l(n+1)}}{\sqrt\pi} \sum_{j\geq (n+1)/2} (2\si)^{2j}
\sqrt{\frac 2 {n+1}} \left(\sqrt{\frac {n+1}{n+3}}\right)^{j-{(n+1)/2}}
,
\end{align*}
with the usual geometric sequence formula giving the statement.
\end{proof}
This section's last method of coping with $j^{-1/2}$
relies upon the chain of equalities\footnote{The idea for the approach comes rather 
naturally if attempting to relate \ref{eq:ctail_stirling} to the Chernoff bound, as addressed
in \autoref{sec:chern}.}
\begin{equation}
\label{eq:chern_funky_sub}
\sum_{j\geq \eta} \frac {x^j}{\sqrt j}
=
\sqrt {\left(\sum_{j\geq \eta} \frac {x^j}{\sqrt j}\right)
\left(\sum_{j\geq \eta} \frac {x^j}{\sqrt j}\right)}
=
\sqrt{
\sum_{k\geq 2\eta}x^k\sum_{j=\eta}^{k-\eta} \frac 1 {\sqrt{j(k-j)}}
}.
\end{equation}
\begin{defn}
For any $\eta,k\in \Z^+$ with $2\eta \leq k$, define
\begin{align*}
\psi_\eta(k) 
&= \sum_{j=\eta}^{k-\eta} \frac 1 {\sqrt{j(k-j)}}.
\end{align*}
\end{defn}
As it turns out, $\psi_\eta$ is rather well behaved.
\begin{lem}
\label{lem:psi_properties}
For any $\eta,k\in \Z^+$ with $2\eta \leq k$, $\psi_\eta(k) \leq \psi_\eta(k+1) \leq \pi$.
\end{lem}
\begin{proof}
Note that
\begin{align*}
\psi_\eta(k)=\sum_{j=\eta}^{k-\eta} \frac 1 {\sqrt{j(k-j)}}
&= \sum_{j=\eta}^{k-\eta} \frac {d}{dj} \cos^{-1}(1-\frac {2j}k) \\
&= -\frac 2 k\sum_{\substack{i = 1 - 2j/k \\ j=\eta,\ldots,k-\eta}} \frac {d}{di}\cos^{-1}(i)
= \frac 2 k\sum_{\substack{i = 1 - 2j/k \\ j=\eta,\ldots,k-\eta}} \frac 1 {\sqrt{1-i^2}}.
\end{align*}
As such, $\psi_\eta(k)$ can be interpreted as a Riemann sum lower bounding the function $(1-x^2)^{-1/2}$ on the interval $(-1,+1)$.  Indeed, take $2/k$ to be the width of each 
rectangle, and when $i\leq 0$, take $(1+i^2)^{-1/2}$ to be the height at the right 
endpoint of a rectangle, otherwise when $i > 0$ take it to be the height at the left
endpoint.  Since $d(\cos^{-1}(x))/dx = -(1-x^2)^{-1/2}$, the value of the approximated 
integral is $\pi$, which gives the upper bound.  (To handle the discontinuity, apply
the monotone convergence theorem to $\lim_{n\to\infty} \int_{-1+1/n}^{1-1/n} \frac {dx}{\sqrt{1-x^2}}$.)

For the monotonicity statement, note that the Riemann sums of a convex, decreasing
function are increasing as the width of the subdivisions decreases (a proof of this
fact is in \autoref{sec:riemann_lemma}).  The result follows by applying this
to both halves of the function separately.
\end{proof}
With this machinery in place, the final bounds of this section may be established.
\begin{thm}
When $n$ is odd
and $p < 1/2$,
\begin{equation}
\label{eq:ctail_chernlike}
\frac {(2\si)^{n+1}e^{l(n+1)}}{\sqrt{2\pi(n+1)}}
\leq
\P[B(p,n)\geq n/2]
\leq
\frac {(2\si)^{n+1}}{2}.
\end{equation}
\end{thm}
It is shown in \autoref{sec:chern} that, with an even number of trials, the
upper bound of \ref{eq:ctail_chernlike} is tighter than the Chernoff bound for all $p\in[0,1/2]$.
Furthermore, the ratio of the two approaches 2 as the number of trials grows.

The lower bound of \ref{eq:ctail_chernlike} has the weakness that, as the
number of trials grows, it becomes poor.
On the other hand, it has the
distinction, among all bounds of this section derived from \ref{eq:ctail_stirling}, that it does
not approach 0 as $p\to 1/2$.

\begin{proof}
Again using \ref{eq:ctail_stirling} but now dealing with $j^{-1/2}$ via \ref{eq:chern_funky_sub},
\begin{align*}
\P[B(p,n)\geq n/2]
&\leq
\frac {1-2p}{2\sqrt{\pi}}
\sqrt{
\sum_{k\geq n+1}
(2\si)^{2k}\psi_{(n+1)/2}(k)} \\
&\leq
\frac {1-2p}{2}
\sqrt{
\frac {(2\si)^{2n+2}}{1-4\si^2}},
\end{align*}
where the conclusion used $\psi_\eta(k) \leq \pi$; to finish, substitute $\sqrt{1-4\si^2} = |1-2p|$.
The lower bound proceeds analogously, but invoking 
\autoref{lem:psi_properties} to grant
$2/(n+1) = \psi_{(n+1)/2}(n+1) \leq \psi_{(n+1)/2}(k)$ for all $k\geq n+1$.
\end{proof}

Note that all preceding bounds used the infinite summation form as presented in
\autoref{lem:ctail_infinite}.  For the last pair of bounds in this section, the
finite sum from \autoref{lem:ctail_finite} is used, which predictably leads to
a much different
bound.

\begin{thm}
When $n$ is odd
and $p < 1/2$,
\begin{align}
\P[B(p,n)\geq n/2]
&\geq
p - \frac{2\si^2e^{u(n-1)}\sqrt{(1 - (2\si)^{2n-4})\psi_1(n-1)}} {\sqrt \pi}
\label{eq:weird_finite_lb},\\
\P[B(p,n)\geq n/2]
&\leq
p - \frac{2\si^2e^{l(2)}\sqrt{1 - (2\si)^{n-1}}} {\sqrt \pi}.
\label{eq:weird_finite_ub}
\end{align}
\end{thm}
Both \ref{eq:weird_finite_lb} and \ref{eq:weird_finite_ub} become exact as $p\to 1/2$, but
are otherwise inaccurate.
\begin{proof}
This proof does not differ greatly from the others in this section, with the exception of starting 
from \autoref{lem:ctail_finite}.  For the lower bound, the key inequalities are
\[
-\sum_{j=1}^{(n-1)/2}\frac {(2\si)^{2j}}{\sqrt j}
\geq -\sqrt{
\sum_{k=2}^{n-1}(2\si)^{2k}\psi_1(k)}
\geq -\sqrt{
\frac 
{((2\si)^4 - (2\si)^{2n})\psi_1(n-1)}{1-4\si^2}},
\]
which makes use of the monotonicity of $\psi_\eta$.  The upper bound is 
similar, but using the fact that $\psi_1\geq 1$.
\end{proof}


\section{Relationship to the Chernoff Bound}
\label{sec:chern}
Since $[X \geq a] = [e^{tX} \geq e^{ta}]$ for all $t>0$, it follows by Markov's inequality
that
\[
\P[X \geq a] \leq \inf_{t>0} \frac{\E(e^{tX})}{e^{ta}};
\]
this is a form of the Chernoff bound (see (3.6) in \citet{chernoff}).
Applying this to the central tail, 
when $m$ is even,
yields
\[
\P[B(p,m) \geq m/2] \leq (2\si)^m.
\]
It is no coincidence this bears a striking resemblance to the upper bound in
\ref{eq:ctail_chernlike}; that bound was derived with the intent of proving the Chernoff
bound.  In fact, adjusting \ref{eq:ctail_chernlike} to 
even $m$
as per \autoref{rem:ctail_constraints},
\[
\P[B(p,m)\geq m/2] 
\leq \frac {(2\si)^m} 2 + \frac 1 2 \binom{m}{m/2}\si^m
\leq (2\si)^m\left(\frac 1 2 + \frac 1 {\sqrt{2\pi m}}\right).
\]
This serves to not only prove the Chernoff bound (for this case), it also states that the multiplicative
error of the Chernoff bound is at least $(1/2 + (2\pi m)^{-1/2})^{-1}$.

When $n$ is odd,
the Chernoff bound is the slightly uglier expression
\begin{equation}
\P[B(p,n)\geq (n+1)/2] \leq (2\sigma)^n 
\sqrt{ \frac {p}{1-p} } \left(\frac {n^2}{n^2-1}\right)^{n/2}\sqrt{\frac{n-1}{n+1}}.
\end{equation}
The ratio of this expression to the upper bound in \ref{eq:ctail_chernlike} is
\begin{equation}
\label{eq:chern_n_odd_p_over_quarter}
\sqrt{\frac 1 {1-p}} \left(\frac {n^2}{n^2-1}\right)^{n/2}\sqrt{\frac{n-1}{n+1}},
\end{equation}
which approaches $2$ as $n\to\infty$.  
It does not, however, exceed 1 for all $p$ and all $n$.  On the other hand, dividing 
the bound by the upper bound in \ref{eq:elem_ub} yields the ratio
\begin{equation}
\label{eq:chern_n_odd_p_under_quarter}
\frac {(1-2p)\sqrt{\pi(n-1)}}{(1-p)\sqrt{2}} \left(\frac {n^2}{n^2-1}\right)^{n/2}.
\end{equation}
The ratio in \ref{eq:chern_n_odd_p_over_quarter} holds when $p\geq 1/4$, whereas the latter ratio 
in \ref{eq:chern_n_odd_p_under_quarter} exceeds 1
when $p\leq 1/4$.  Thus, combining the two via a $\min$ yields a better bound.  As will be discussed in \autoref{sec:allbounds}, a number of the bounds, when paired via $\min$
or $\max$, form extremely good bounds.

Lastly, note that the Chernoff bound was chosen because, for an even number of trials, it is
it is tighter than the corresponding Hoeffding and Bernstein bounds.
(For a proof, see 
\autoref{sec:pf_chern_cmp}.)

\section{Bounding with the Standard Normal}
\label{sec:contbounds}
The preceding bounds all aimed for a closed-form approximation for either the
infinite or finite summation in \autoref{thm:ctail}. 
In this section, however, the strategy is to replace the infinite summation with an integral.  
As usual,
let $\Phi$ and $\phi$ be the distribution function and density of the standard normal.

\begin{thm}
\label{thm:cont_bounds}
Let $n$ odd
and $p\in (0,1/2)$ be given, and set
\begin{align*}
\Upsilon 
&= \frac{2\left(1-\Phi(\sqrt{-(n+1)\ln(4\si^2)})\right)}{\sqrt{-\ln(4\si^2)}},
&\Delta 
&= \frac {(2\si)^{n+1}}{\sqrt {2\pi(n+1) }};
\end{align*}
then
\begin{equation}
\label{eq:cont_bounds}
(1/2-p)e^{l(n+1)}(\Upsilon + \Delta)
\leq
\P[B(p,n)\geq n/2]
\leq
(1/2-p)(\Upsilon + \Delta(1+R)),
\end{equation}
where 
\[
R 
\leq \min\left\{
1,
\frac 1 4 \left( \frac 1 {n+1} - \ln (4\si^2) \right)
\right\}.
\]
\end{thm}
The most important property is that the expressions for the upper and lower bounds are nearly 
the same, providing for easy comparison.  Concretely, the additive error of either can be bounded with
their difference
\begin{equation*}
(1/2 - p)(\Upsilon + \Delta(1+R)) - (1/2 - p)e^{l(n+1)}(\Upsilon + \Delta).
\end{equation*}
Using $e^{l(n+1)} \geq e^{-1/3(n+1)} \geq 1 - 1/3(n+1)$ and $R\leq 1$
yields
\begin{equation}
\label{eq:contb_2}
\frac 1 {3(n+1)} \left(\frac 1 2 - p\right) \Upsilon + \left(1 + \frac 1 {3(n+1)}\right)\left(\frac 1 2 - p\right)\Delta.
\end{equation}
Since  $(1/2 - p) \Upsilon$ is increasing (along $p\in[0,1/2]$), it may be replaced with its
limiting value $1/2$.
Substituting the maximizing value for $p$ into the right summand of \ref{eq:contb_2} and
simplifying, the error is thus upper bounded by $2/5(n+1)$.
Although coarse, this error bound 
provides some explanation of the accuracy of \ref{eq:cont_bounds}, evidenced in
\autoref{fig:smallplot} and \autoref{fig:plots}.

The plots use the minimum of the two choices for $R$ (for every $p$).
As is discussed in \autoref{sec:slud}, it is possible to prove Slud's bound
by relaxing these bounds, and when a choice for $R$ must be made, 1 suffices.  Lastly note that the more
complicated bound on $R$, though usually better, is worse for small values of $p$.  For instance, a sufficient
condition for the complicated bound to be better (for any $n\geq 1$) is $p \geq 0.0077$.

To prove the theorem, first note that a change of variable suffices to remove $j^{-1/2}$ from the
integral.  (Recall that this term was the primary difficulty in \autoref{sec:elembounds}.)

\begin{lem}
\label{lem:int_bound_helper}
When $\si\in(0,1/2)$,
\[
\int_{(n+1)/2}^\infty \frac {(2\si)^{2j}}{\sqrt{\pi j}}dj
=
\frac{2\left(1-\Phi(\sqrt{-(n+1)\ln(4\si^2)})\right)}{\sqrt{-\ln(4\si^2)}}.
\]
\end{lem}
\begin{proof}
To start,
\[
\int_{(n+1)/2}^\infty \frac {(2\si)^{2j}}{\sqrt{\pi j}}dj
= \int_{(n+1)/2}^\infty \frac {e^{j\ln(4\si^2)}}{\sqrt{\pi j}}dj.
\]
Applying the map $j\mapsto -j^2 / (2\ln(4\si^2))$ yields
\[
\int_{\sqrt{-(n+1)\ln(4\si^2)}}^\infty \frac {e^{-j^2/2}}{\sqrt{-\pi j^2 / (2\ln(4\si^2))}}\left(\frac{-2j}{2\ln(4\si^2)}\right)dj,
\]
which gives the statement after some algebra.
\end{proof}

\begin{proof}[Proof of \autoref{thm:cont_bounds}]
By the first order Euler-Maclaurin summation formula,
\[
\sum_{j\geq(n+1)/2}\!\!\!\!\! \frac {(2\si)^{2j}}{\sqrt{\pi j}}
= \int_{(n+1)/2}^\infty\!\!\!\!\frac {(2\si)^{2j}}{\sqrt{ \pi j}} dj
+\frac {(2\si)^{n+1}}{\sqrt{2\pi(n+1)}}
+\int_{(n+1)/2}^\infty\!\!\!\!\!\!\!\!\!\! (\{x\}-1/2)\left(\frac d {dx}\frac {(2\si)^{2x}}{\sqrt{\pi x}}\right)dx, 
\]
where $\{x\}$ is the fractional part of $x$.  (To turn this expression into lower and upper bounds for $\P[B(p,n)\geq n/2]$,
scale by 
$(1-2p)e^{l(n+1)}$
and $1-2p$, which yields both sides of \ref{eq:ctail_stirling}).  By \autoref{lem:int_bound_helper}, the first term is
$\Upsilon$, and the second is $\Delta$, so only the last term requires attention.  To start, observe that
\begin{align*}
\frac {d}{dx} \left(\frac{(2\si)^{2x}}{\sqrt{\pi x}}\right)
&= \frac{(2\si)^{2x}}{\sqrt{\pi x}} \left(\ln(4\si^2) - \frac 1 {2x}   \right) =:f(x), \\
\frac {d^2}{dx^2} \left(\frac{(2\si)^{2x}}{\sqrt{\pi x}}\right)
&= \frac{(2\si)^{2x}}{\sqrt{\pi x}}\left(
\ln^2(4\si^2) - \frac {\ln(4\si^2)}{x} + \frac {3}{4x^2}
\right).
\end{align*}
Since $\si\in(0,1/2)$, $f\leq 0$ and $f'\geq 0$, and hence the
integral is nonnegative (once again, this follows from exercise 9.16 of \citet{GKP}), thus establishing the
lower bound. Next
\begin{align*}
\int_{(n+1)/2}^\infty (\{x\}-1/2)\left(\frac d {dx}\frac {(2\si)^{2x}}{\sqrt{\pi x}}\right)dx
&\leq
-\frac 1 2\int_{(n+1)/2}^\infty \left(\frac d {dx}\frac {(2\si)^{2x}}{\sqrt{\pi x}}\right)dx
= \Delta,
\end{align*}
establishing the first upper bound on $R$.  (The same bound may be derived by
starting with the naive integral bound instead of Euler-Maclaurin.)  For the
second upper bound, write
\begin{align*}
\int_{(n+1)/2}^\infty(\{x\}-1/2) f(x) dx
&=\sum_{j\geq (n+1)/2} \int_0^1(x-1/2) f(j+x) dx\\
&\hspace{-1.25in}=\sum_{j\geq (n+1)/2}\left( \int_0^{1/2}(x-1/2) f(j+x) dx
+\int_{1/2}^1(x-1/2) f(j+x) dx
\right).
\end{align*}
Since $f(x)$ is negative and monotonic increasing, it follows that
\begin{align*}
\int_0^{1/2}(x-1/2) f(j+x) dx
&\leq \inf_{x\in[0,1/2]} f(j+x) \int_0^{1/2}(x-1/2) dx = -f(j)/8,\\
\int_{1/2}^1(x-1/2) f(j+x) dx
&\leq \sup_{x\in[1/2,1]} f(j+x) \int_{1/2}^1(x-1/2) dx = f(j+1)/8.
\end{align*}
As such, the sum telescopes, establishing the other bound on $R$.
\end{proof}
\section{Relationship to Slud's Bound}
\label{sec:slud}
Slud's bound is the standard tool for lower bounding binomial tails.

\begin{thm}[\citet{slud}] 
Let $n,k$ be nonnegative integers with $k\leq n$, and $p\in[0,1]$.  When either
(a) $p\leq 1/4$ and $np\leq k \leq n$, or (b) $np \leq k \leq n(1-p)$, then
\begin{equation}
\P[B(p,n) \geq k] \geq 1- \Phi\left(\frac {k-np}{\si\sqrt{n}}\right).
\end{equation}
\end{thm}
\begin{rem}
Many presentations omit sufficient condition (a).  Many also omit the integrality of $k$, which can
be seen as necessary with the example $n=1$ and $p=k\in(0,0.5)$; in this case, Slud's lower
bound is exactly $0.5 > p = \P[B(p,1)\geq k]$.
\end{rem}

It is possible to start with the bounds in \autoref{thm:cont_bounds}, and
apply a battery of elementary inequalities (mostly tangent and secant approximations to
the relevant functions) to weaken the inequalities into Slud's bound.  The proof
is quite tedious, and thus deferred to 
\autoref{sec:pf_slud},
however a few points are worthy 
of mention.

By the method of proof, it is immediate that the bounds of \autoref{thm:cont_bounds} are tighter
than Slud's inequality.  In one case, care is even made to maintain a small separation, however quantifying
the exact gap is hard.

Perhaps most importantly, the proof was able to extend the sufficient conditions for Slud's inequality.
For $n$ odd and $m$ even, Slud's bound (for central tails) requires $p\leq1/2-1/2n$ and $p\leq 1/2$, respectively.
The new proof, however, holds for 
\begin{align*}
p &\leq \frac 1 2 + \frac 1 2 \left(1-\frac{4(\sqrt{n(n+1)} - 1)}{4n+2}\right)^{1/2},
&p &\leq \frac 1 2 + \frac 1 6(e^{l(m)}/m)^{1/3}.
\end{align*}
This is nice since Slud's bound is of a significantly simpler (and more interpretable) form than the bounds
of \autoref{thm:cont_bounds}.

Empirical evidence seems to suggest that Slud's bound does not hold for all $p$, and in fact, as $n\to\infty$,
the maximal permissible $p$ shrinks to $1/2$.  Also, the following appears to be true.

\begin{conj}
When $p\in(0,1/2)$ and 
$m$ is even, 
\[
\P[B(p,m) \geq m/2] \geq 1 - \Phi\left(\frac{(1/2-p)\sqrt{m}}{\si}\right) + \frac 1 2\binom{m}{m/2}\si^{m}.
\]
\end{conj}
Unfortunately, the bounds of \autoref{thm:cont_bounds} are not sufficiently tight to establish
this; perhaps if the term $e^{l(n+1)}$ were handled better.



\section{Summary of Bounds}
\label{sec:allbounds}
\begin{figure}
\includegraphics[width=\textwidth]{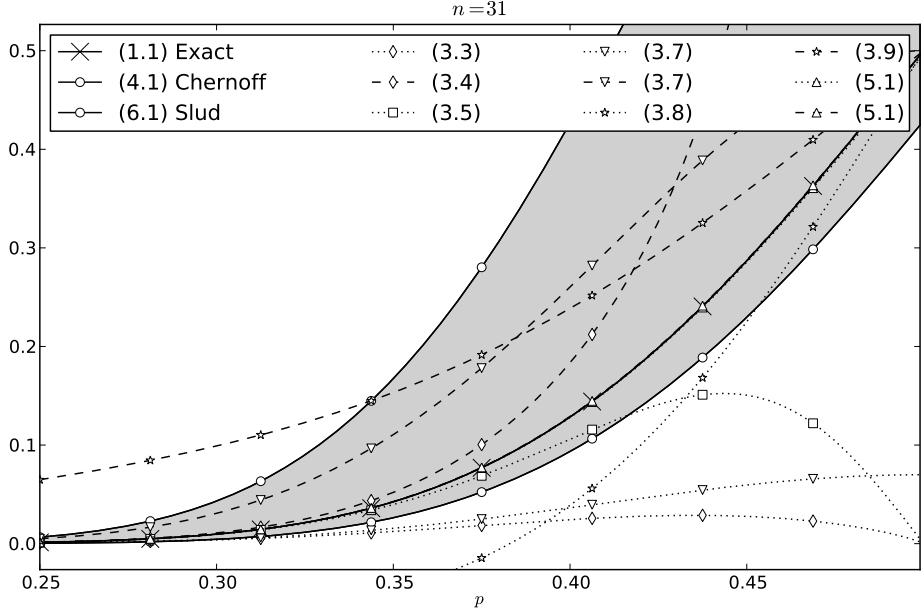}
\caption{
Bounds are identified by equation number, lower and upper bounds derived in this paper being
respectively dotted and dashed.  Shading identifies
the region between the Chernoff and Slud bounds.
The exact (central) binomial tail is
marked with x's.
}
\label{fig:plots}
\end{figure}

\autoref{fig:plots} contains plots of all bounds.
The bounds of \ref{eq:cont_bounds} are almost exact, which is in agreement with 
their error bound.
The simple upper bound in \ref{eq:ctail_chernlike} is also generally good.  On the other hand,
most of the other bounds vary performance quite widely with $p$.
This suggests use of pairs of bounds in tandem;
for instance, as was mentioned in \autoref{sec:chern}, the minimum
of the upper bounds in \ref{eq:ctail_chernlike} and \ref{eq:elem_ub}
suffices to prove Chernoff's bound.
Similarly, the $\max$ of the lower bounds in \ref{eq:funky_exp_lower_bound} and 
\ref{eq:weird_finite_lb} would work well.



Slud's bound fares well, being the best lower bound with the exception of the
bounds in \ref{eq:cont_bounds}, and \ref{eq:funky_exp_lower_bound} for certain values of $p$.
In contrast, many upper bounds outperform the Chernoff bound.

\section{General Binomial Tails via Random Walks}
\label{sec:gtail}
By following steps analogous to those of \autoref{sec:basics}, \autoref{thm:ctail}
readily generalizes. (For details, see 
\autoref{sec:pf_gtail}.)

\begin{thm}
\label{thm:gtail}
When $n$ is odd,
$k\in \Z \cap [0,n/2)$, and $p \in (0,1/2)$,
\begin{align}
\P[B(p,n)\geq \frac n 2 + k] 
&= p^{2k+1}
- \left(\frac {p}{1-p}\right)^k\sum_{j=k+1}^{(n-1)/2} 
\left(\frac 1 2 -p - \frac{k}{2j}\right)\binom{2j}{j+k} \si^{2j}\label{eq:gtail_a}\\
&=\left(\frac {p}{1-p}\right)^k \sum_{j\geq (n+1)/2}\left(\frac 1 2 -p - \frac k {2j}\right)
\binom {2j}{j+k} \si^{2j}.\label{eq:gtail_b}
\end{align}
Furthermore
\begin{align}
\P[B(1/2,n)\geq \frac n 2 + k] 
&= 2^{-2k-1}
+ \sum_{j=k+1}^{(n-1)/2} 
\left(\frac{k}{2j}\right)\binom{2j}{j+k} 4^{-j}\label{eq:gtail_c}\\
&= \frac 1 2
- \sum_{j\geq (n+1)/2}\left(\frac k {2j}\right)
\binom {2j}{j+k} 4^{-j}.\label{eq:gtail_d}
\end{align}
\end{thm}

\section*{Acknowledgements}
The author thanks Prof. Patrick Fitzsimmons for his gracious advice and support; the author
would also like to thank
Daniel Hsu, Jeroen Rombouts, and Michael Wasson for stimulating conversations.

\bibliographystyle{plainnat}

\newcommand{\thmcontbounds}{\autoref{thm:cont_bounds}}
\newcommand{\remctailconstraints}{\autoref{rem:ctail_constraints}} 
\newcommand{\eqcentralstirling}{\ref{eq:central_stirling}} 
\newcommand{\lemctailfinite}{\autoref{lem:ctail_finite}} 

\appendix
\section{Proof of Theorem~\ref{thm:gtail}} 
\label{sec:pf_gtail}
The first task is to generalize \autoref{lem:ctail_step}.
\begin{lem} \label{lem:gtail_step}
When $n$ is odd and $k \in \Z\cap[-n/2,n/2]$,
\begin{align*}
&\P\left[B(p,n+2) \geq \frac{n+2} 2 + k\right] - \P\left[B(p,n) \geq \frac n 2 + k\right] \\
&\quad= \left(p - \frac  1 2+ \frac k {n+1}\right)\binom{n+1}{(n+1)/2 + k}
p^{(n+1)/2 + k} (1-p)^{(n+1)/2 - k}.
\end{align*}
\end{lem}
\begin{proof}
As in the proof of
\autoref{lem:ctail_step},
 the change in probability mass comes entirely from
the random walks which after $n$ steps are at coordinates $(n-1)/2 + k$ or $(n+1)/2 +k$.  Thus,
as before, the mass gained minus the mass lost is
\begin{align*}
&p^2\binom {n}{(n-1)/2 + k} p^{(n-1)/2 + k} (1-p)^{(n+1)/2 -k}\\
&\quad- (1-p)^2\binom {n}{(n+1)/2 + k} p^{(n+1)/2 + k} (1-p)^{(n-1)/2 -k} \\
&=
p^{(n+1)/2 + k} (1-p)^{(n+1)/2 - k}\left[
p\binom {n}{(n-1)/2 + k}
- (1-p) \binom {n}{(n+1)/2 + k}
\right] \\
&=
\binom{n+1}{(n+1)/2+k}p^{(n+1)/2 + k} (1-p)^{(n+1)/2 - k} \\
&\quad\cdot\left[
p\left(\frac{(n+1)/2+k}{n+1}\right)
- (1-p) \left(\frac{(n+1)/2-k}{n+1}\right)
\right] \\
&=
\binom{n+1}{(n+1)/2+k}p^{(n+1)/2 + k} (1-p)^{(n+1)/2 - k}\left[
p- \frac 1 2 + \frac k {n+1}
\right].\qedhere
\end{align*}
\end{proof}
Substituting $k=0$ in the above yields both the statement (and proof) of Lemma 2.1.  
\begin{proof}[Proof of \autoref{thm:gtail}]
Since the case $k=0$ is Theorem 1.1, take $k\neq 0$.
To show \ref{eq:gtail_a}, it suffices to rewrite the tail as a telescoping series, and invoke
\autoref{lem:gtail_step}.  Note that this derivation holds for any $p\in(0,1)$.
\begin{align*}
&\P[B(p,n)\geq n/2 + k] \\
&= \P[B(p,2k+1) \geq (2k+1)/2 + k] \\
&\quad + \sum_{j=1}^{(n-1)/2 - k} \Big(
\P[B(p,2k+1+2j) \geq (2k+1)/2 + k + j]\\
&\quad\quad\quad\quad- \P[B(p,2k-1+2j) \geq (2k-1)/2 + k + j]
\Big) \\
&= p^{2k+1}
+ \sum_{j=1}^{(n-1)/2 - k} 
\left(p- \frac 1 2 + \frac{k}{2j+2k}\right)\binom{2j+2k}{j+2k} p^{2k+j}(1-p)^{j}\\
&= p^{2k+1}
+ \left(\frac {p}{1-p}\right)^k\;\;\sum_{j=k+1}^{(n-1)/2} 
\left(p- \frac 1 2 + \frac{k}{2j}\right)\binom{2j}{j+k} \si^{2j}.
\end{align*}
From here, \ref{eq:gtail_c} may be obtained by substituting $p=1/2$.
Note that, with the exception of the last step, $p$ can happily take on values in $\{0,1\}$,
and the expression in the penultimate yields probabilities of $0,1$, respectively.

Next, using (2.5.16) from \citet{gfology},
\begin{align}
\left(\frac {p}{1-p}\right)^k  \;\;
\sum_{j\geq k}\left(\frac k {2j}\right) \binom {2j}{j+k} \si^{2j}
&= 
\frac {p^{2k}}{2}
\sum_{j\geq 0} \frac {2k}{2j+2k}\binom {2j+2k}{j} \si^{2j} \notag\\
&= 
\frac {p^{2k}}{2}
\left(
\frac {1 - \sqrt{1-4\si^2}}{2\si^2}
\right)^{2k} \notag\\
&= 
\frac {p^{2k}}{2}
\left(
\frac {1 - |2p-1|}{2\si^2}
\right)^{2k}.\label{eq:gtail_b_1}
\end{align}
When $p=1/2$, this expression is $1/2$; thus, starting from \ref{eq:gtail_c},
\begin{align*}
&2^{-2k-1}
+ \sum_{j=k+1}^{(n-1)/2} 
\left(\frac{k}{2j}\right)\binom{2j}{j+k} 4^{-j}\\
&=2^{-2k-1}
+ \sum_{j\geq k}
\left(\frac{k}{2j}\right)\binom{2j}{j+k} 4^{-j}\\
&\quad-
\left(\frac{k}{2k}\right)\binom{2k}{k+k} 4^{-k}
- \sum_{j\geq (n+1)/2}
\left(\frac{k}{2j}\right)\binom{2j}{j+k} 4^{-j} \\
&= \frac 1 2
- \sum_{j\geq (n+1)/2}\left(\frac k {2j}\right)
\binom {2j}{j+k} 4^{-j},
\end{align*}
which is exactly \ref{eq:gtail_d}.  Finally, to handle \ref{eq:gtail_b}, first use
(2.5.15) from \citet{gfology} to obtain
\begin{align}
\left(\frac {p}{1-p}\right)^k  \;\;
\sum_{j\geq k}\left(p - \frac 1 2\right) \binom {2j}{j+k} \si^{2j}
&= 
p^{2k} \left(p - \frac 1 2\right) 
\sum_{j\geq 0} \binom {2j+2k}{j} \si^{2j}\notag\\
&= 
p^{2k} \left(p - \frac 1 2\right) 
\frac {1}{\sqrt{1-4\si^2}}
\left(
\frac {1 - \sqrt{1-4\si^2}}{2\si^2}
\right)^{2k}\notag\\
&=
\left(\frac {p^{2k}}{2}\right)
\frac {2p-1}{|2p-1|}
\left(
\frac {1 - |2p-1|}{2\si^2}
\right)^{2k}.\label{eq:gtail_b_2}
\end{align}
Combining \ref{eq:gtail_b_1} and \ref{eq:gtail_b_2}, it follows that
\begin{align*}
\left(\frac {p}{1-p}\right)^k \;\;\sum_{j\geq k}\left(p - \frac 1 2 + \frac k {2j}\right)
\binom {2j}{j+k} \si^{2j}
&= 
\frac {p^{2k}}{2}
\left(
\frac {1 - |2p-1|}{2\si^2}
\right)^{2k}
\left(
\frac {2p-1}{|2p-1|}
+1
\right),
\end{align*}
which further reduces to just $0$ when $p < 1/2$.
Thus, from \ref{eq:gtail_a},
\begin{align*}
\P[B(p,n)\geq n/2 + k]
&= p^{2k+1}
- \left(\frac {p}{1-p}\right)^k\;\;\sum_{j=k+1}^{(n-1)/2} 
\left(\frac 1 2 -p - \frac{k}{2j}\right)\binom{2j}{j+k} \si^{2j}\\
&= p^{2k+1}
- \left(\frac {p}{1-p}\right)^k\;\;\sum_{j\geq  k}
\left(\frac 1 2 -p - \frac{k}{2j}\right)\binom{2j}{j+k} \si^{2j}\\
&\quad+ \left(\frac {p}{1-p}\right)^k
\left(\frac 1 2 -p - \frac{k}{2k}\right)\binom{2k}{k+k} \si^{2k}\\
&\quad+ \left(\frac {p}{1-p}\right)^k\sum_{j\geq (n+1)/2}
\left(\frac 1 2 -p - \frac{k}{2j}\right)\binom{2j}{j+k} \si^{2j}\\
&= p^{2k+1} - 0 - p^{2k+1}\\
&\quad + \left(\frac {p}{1-p}\right)^k\sum_{j\geq (n+1)/2}
\left(\frac 1 2 -p - \frac{k}{2j}\right)\binom{2j}{j+k} \si^{2j},
\end{align*}
which is \ref{eq:gtail_b}.
\end{proof}

\section{Proof of (Central) Slud's Bound}
\label{sec:pf_slud}

\begin{thm}
When $n$ is odd and
$0\leq p \leq \frac 1 2 + \frac 1 2 \left(1-\frac{4(\sqrt{n(n+1)} - 1)}{4n+2}\right)^{1/2}$,
\begin{equation}
\label{eq:slud_odd_bound}
\P[B(p,n) \geq n/2] \geq 1 - \Phi\left(\frac {(n+1)/2 - np}{\si\sqrt{n}}\right).
\end{equation}
When $m$ is even and
$0 \leq p \leq \frac 1 2 + \frac 1 6(e^{l(m)}/m)^{1/3}$,
\begin{equation}
\label{eq:slud_even_bound}
\P[B(p,m/2)\geq m/2] \geq 1- \Phi\left(\frac {(1/2 - p)\sqrt{m}}{\si}\right).
\end{equation}
\end{thm}
When $p\in\{0,1/2\}$, the statements are immediate, and thus disregarded.  The
proof is split into four parts, each reducing to the
bounds in \thmcontbounds.

\begin{proof}[Proof of \ref{eq:slud_odd_bound} when $0<p<1/2$]
Set $\alpha = \sqrt{-(n+1)\ln(4\si^2)}$ and $\beta = ((n+1)/2 - np)/(\si\sqrt{n})$; note that
$\alpha < \beta$.  Using the lower bound in \thmcontbounds. the statement is
implied by
\[
\frac {(1-2p)e^{l(n+1)}(1-\Phi(\alpha))}{\sqrt{-\ln(4\si^2)}}
+ \frac {(1-2p)e^{l(n+1)}(2\si)^{n+1}}{\sqrt {2\pi(n+1) }}
\geq 1 - \Phi(\beta).
\]
Dropping the second term and using $\sqrt{-\ln (4\si^2)} \leq
\sqrt{(4\si^2)^{-1} -1} = 2\si/(1-2p)$, this is a consequence of
\[
2\si e^{l(n+1)}(1-\Phi(\alpha))\geq 1- \Phi(\beta).
\]
$0\leq \alpha \leq \beta$ so $\Phi(\beta)-\Phi(\alpha) \geq \phi(\beta)(\beta-\alpha)$, and
$1-\Phi(\beta) < \phi(\beta)/\beta)$ 
(cf. (1.8) from chapter 7 of \citet{feller}),
 so this in turn is
implied by either of
\begin{align*}
2\si e^{l(n+1)}\phi(\beta)(\beta-\alpha) 
&\geq (1-2\si e^{l(n+1)})\frac {\phi(\beta)}{\beta}
&\iff\\
2\si e^{l(n+1)}(\beta(\beta -\alpha)+1) 
&\geq 1.
\end{align*}
Establishing the latter inequality is \autoref{lem:shitty_derivative_slud_shit}.
\end{proof}

\begin{proof}[Proof of \ref{eq:slud_odd_bound} when 
$\frac 1 2<p<\frac 1 2 + \frac 1 2 \left(1-\frac{4(\sqrt{n(n+1)} - 1)}{4n+2}\right)^{1/2}$]
Set 
$\beta = (np - (n+1)/2)/(\si\sqrt{n})$
and 
$\alpha = \sqrt{-(n+1)\ln(4\si^2)}$.
By
\remctailconstraints\ 
 and symmetry of $\phi$, the theorem statement may be rewritten as
\[
\P[B(1-p,n)\geq n/2] 
\leq 1 - \Phi(\beta).
\]
Discard
the case that $\beta<0$, since then $1-\Phi(\beta)>1/2$ whereas $\P[B(1-p,n)\geq n/2]\leq 1/2$ 
by \lemctailfinite.  Since $2p-1 \leq \sqrt{-\ln(4\si^2)}$ this in turn is a consequence
of
\begin{equation}
\label{eq:slud_odd_2_1}
1-\Phi(\alpha) + \frac{(1-2p)(2\si)^{n+1}}{\sqrt{2\pi(n+1)}} 
\leq 1 -\Phi(\beta).
\end{equation}
As per \autoref{lem:other_shitty_slud_shitty}, the conditions on $p$ imply
\[
\frac {2p-1}{\sqrt{n+1}} \leq \alpha - \beta,
\]
meaning $\alpha\geq \beta$.  Thus scaling both  sides by $\phi(\alpha)$ yields
\[
\frac {(2p-1)(2\si)^{n+1}}{\sqrt{2\pi(n+1)}} \leq (\alpha - \beta)\phi(\alpha),
\]
which implies \ref{eq:slud_odd_2_1} (since an integral can be lower bounded by a rectangle).
\end{proof}

\begin{proof}[Proof of \ref{eq:slud_even_bound} when $0<p<1/2$]
Set $\alpha = \sqrt{-m\ln(4\si^2)}$ and $\beta = ((1/2 - p)\sqrt{m})/\si$; again, $\alpha \leq \beta$. 
Invoking \remctailconstraints, \eqcentralstirling, and \thmcontbounds, the
theorem statement is implied by
\[
e^{l(m)}\left(\frac{(1-2p)(1-\Phi(\alpha))}{\sqrt{-\ln(4\si^2)}} + \frac{(3/2-p)(2\si)^m}{\sqrt{2\pi m}}\right) \geq 1-\Phi(\beta).
\]
As in the proof 
when $n$ is odd, use 
$\sqrt{-\ln (4\si^2)} \leq 2\si/(1-2p)$, but simply drop the term $\Phi(\beta)-\Phi(\alpha)$, which
gives the antecedent statement
\begin{equation}
\label{eq:slud_even_1_1}
\frac {e^{l(m)}(3/2-p)(2\si)^m}{\sqrt{2\pi m}} \geq (1 - 2e^{l(m)}\si)(1-\Phi(\beta)),
\end{equation}
which will be established using the fact $1-\Phi(\beta) \leq \min\{\phi(\beta)/\beta, 1/2\}$.
Set $\bar p = \sqrt{3/4 - 1/(2e^{l(m)})}$; 
the statement is established for $p\in[\bar p, 1/2]$ in \autoref{lem:slud_fuckballs},
so take
 $p \in (0,\bar p]$.  Rearranging, the condition on $p$ states
\[
e^{l(m)}(\frac 3 4 - p^2)\geq \frac 1 2,
\]
and $\si \leq 1/2$, so this becomes 
\[
\frac {e^{l(m)}(\frac 3 4 - p^2)}{\si}\geq 1.
\]
Next, $3/4-p^2 = (3/2-p)(1/2-p) + 2\si^2$, so the above can be re-arranged into
\[
\frac {e^{l(m)}(1.5-p)(0.5-p)}{\sigma} \geq 1-2e^{l(m)}\sigma.
\]
To finish, the definition of $\beta$ and scaling by $\phi(\beta)$ implies
\[
\frac {e^{l(m)}(1.5-p)\phi(\beta)}{\sqrt{m}} \geq (1-2e^{l(m)}\sigma)\frac {\phi(\beta)}{\beta},
\]
where substituting $\phi(\beta)\leq \phi(\alpha) = (2\si)^2 / \sqrt{2m}$ 
yields
\ref{eq:slud_even_1_1}.
\end{proof}

\begin{proof}[Proof of \ref{eq:slud_even_bound} when $\frac 1 2<p<1/2 + (e^{l(m)}/m)^{1/3}/6$]
Set $\alpha = \sqrt{-m\ln(4\si^2)}$ and \linebreak$\beta = ((p-1/2)\sqrt{m})/\si$; again 
using \remctailconstraints\ and 
\eqcentralstirling, the theorem statement may be rewritten
\[
\P[B(1-p,m-1) \geq m/2] - \frac {e^{l(m)} (2\si)^{m/2}}{\sqrt{2\pi m}}\leq 
1-\Phi(\beta).
\]
Invoking \thmcontbounds, $2p-1 \leq \sqrt{-\ln(4\si^2)}$, and
$\phi(\alpha) = (2\si)^m/\sqrt{2\pi}$, this is implied by
\[
\frac {(2p-1-e^{l(m)})\phi(\alpha)}{\sqrt{m}} \leq\Phi(\alpha)  - \Phi(\beta),
\]
which by $\Phi(\alpha) - \Phi(\beta) \geq -(\beta-\alpha)\phi(\alpha)$ is a consequence of
\begin{equation}
\label{eq:slud_even_2_1}
\frac {(2p-1-e^{l(m)})}{\sqrt{m}} \leq \alpha - \beta.
\end{equation}
Set $q = 2p-1$; the conditions on $p$ mean $q \leq (e^{l(m)}/m)^{1/3}/3$.
Using algebra and some bounds on $l(m)$, it follows that
\[
(m-1)q^3 + e^{l(m)}q^2 + q -e^{l(m)} \leq 0,
\]
which can be rearranged into
\[
e^{l(m)} - q \geq m\left(\frac {q}{1-q^2} -q \right).
\]
As before, $q=2p-1 \leq \sqrt{-\ln(4\si^2)}$, and since $1-q^2 = 2p(2-p)\leq 2\si$,
this implies
\[
e^{l(m)} - (2p-1) \geq m\left(\frac {2p-1}{2\si} -\sqrt{-\ln(4\si^2)} \right),
\]
which gives \ref{eq:slud_even_2_1}.
\end{proof}

\section{Supporting Lemmas}
\begin{lem}
\label{lem:shitty_derivative_slud_shit}
When $p\in(0,1/2)$ and $n$ odd,
\[
2\si e^{l(n+1)}
\left(
\left(
\frac {(n+1)/2 - np}{\si\sqrt{n}}
\right)
\left(
\frac {(n+1)/2 - np}{\si\sqrt{n}}
-\sqrt{-(n+1)\ln(4\si^2)}
\right)
+1
\right)
\geq 1.
\]
\end{lem}
\begin{proof}
Set $r = p(1-p) = \si^2$, and lower bound the left hand side according to
\begin{align}
&
2\sqrt{r} e^{l(n+1)} \left[
\frac {n\sqrt{1-4r} + 1}{2\sqrt{nr}}
\left(\frac {n\sqrt{1-4r} + 1}{2\sqrt{nr}} - \sqrt{-(n+1)\ln(4r)}\right) + 1
\right] \notag\\
&\geq
2\sqrt{r} e^{l(n+1)} \left[
\frac {n\sqrt{1-4r} + 1}{2\sqrt{nr}}
\left(\frac {n\sqrt{1-4r} + 1}{2\sqrt{nr}} - \frac {\sqrt{(n+1)(1-4r)}}{\sqrt{4r}}\right) + 1
\right] \notag\\
&=
2\sqrt{r} e^{l(n+1)} \Bigg[
\frac {n^2(1-4r)+2n\sqrt{1-4r} + 1}{4nr}\notag\\
&\quad-\frac{n(1-4r)\sqrt{n+1} +\sqrt{(n+1)(1-4r)}}{4r\sqrt n} +1
\Bigg] \notag\\
&=
e^{l(n+1)} \Bigg[
\frac {1-4r}{\sqrt r} \left(\frac n 2 -\frac {\sqrt{n(n+1)}} 2 \right)\label{eq:slud_proof_step_4}\\
&\quad+\sqrt{\frac {1-4r}{r}} \left(1 - \frac 1 2 \sqrt{\frac {n+1}{n}} \right)
+\frac 1 {\sqrt r}\left(\frac 1 {2n}\right)
+\sqrt r(2)\Bigg].\notag
\end{align}
The goal will be to show that this final expression is lower bounded by $1$.  In particular,
its derivative with respect to $r$ is
\begin{align*}
r^{-1/2}e^{l(n+1)}\Bigg[
&-\frac {1+4r}{2r} \left(\frac n 2 -\frac {\sqrt{n(n+1)}} 2 \right)\\
&-\frac {1}{2r\sqrt{1-4r}} \left(1 - \frac 1 2 \sqrt{\frac {n+1}{n}} \right)
-\frac 1 {4nr}
+\frac 1 2 
\Bigg];
\end{align*}
by showing this is negative, it will suffice to show that \ref{eq:slud_proof_step_4} holds
at $p=1/2$.

First consider the case that $n\in\{1,2,3,4\}$ (and note that $r\in(0,1/4)$, since $p\in(0,1/2)$).
Each case can be checked manually by plugging in the given $n$, and using $(1-4r)^{-1/2} \geq 1+2r$ (the linear approximation at $r=0$).

So take $n\geq 5$, and
it will be shown that 
\[
\star = -\frac {1+4r}{2r} \left(\frac n 2 -\frac {\sqrt{n(n+1)}} 2 \right)
-\frac {1}{2r\sqrt{1-4r}} \left(1 - \frac 1 2 \sqrt{\frac {n+1}{n}} \right)
\leq-1/2.
\]
Since $n\geq 5$, using some calculus, it holds that
\begin{align*}
\frac n 2 -\frac {\sqrt{n(n+1)}} 2 &\geq -\frac 1 4,\\
1 - \frac 1 2 \sqrt{\frac {n+1}{n}} &\geq 0.45.
\end{align*}
Therefore
\[
\star \leq -\frac 1 {2r} \left(
-\frac {1+4r}{4} + \frac {0.45}{\sqrt{1-4r}}
\right).
\]
Consider the case that $r\in[0,1/8]$, and replace $(1-4r)^{-1/2} \geq 1+2r$ (again, 
the linear approximation at $r=0$):
\begin{align*}
-\frac 1 {2r} \left(
-\frac {1+4r}{4} + \frac {0.45}{\sqrt{1-4r}}
\right)
&\leq 
-\frac 1 {2r} \left(
-\frac {1+4r}{4} + (0.45)(1+2r)
\right) \\
&=
-\frac 1 {2r} \left(
\frac 2 {10} - \frac r {10}
\right) \\
&=
\frac 1 {20} - \frac 1 {10r} 
\leq 1 {20} - \frac 8 {10} < 1/2.
\end{align*}
Now take $r\in[1/8,1/4]$, and use the linear approximation $(1-4r)^{-1/2} \geq \sqrt{2}(1 + 4r)$:
\begin{align*}
-\frac 1 {2r} \left(
-\frac {1+4r}{4} + \frac {0.45}{\sqrt{1-4r}}
\right) < -1/2.
\end{align*}
Since $\star$ is decreasing, the minimum is obtained when $r=1/4$.  
Plugging $r=1/4$ into \ref{eq:slud_proof_step_4}, 
\begin{align*}
\ref{eq:slud_proof_step_4} 
&\geq e^{l(n+1)} (1+n^{-1}).
\end{align*}
This exceeds 1 for $n=1$; for $n\geq 3$, since 
$1+x \geq \exp(x-x^2)$ when $x\in[0,1/2]$, it follows that
\begin{align*}
e^{l(n+1)} (1+n^{-1})
&\geq \exp(l(n+1)+(n^{-1}-n^{-2}))\\
&\geq \exp(-(3n+1)^{-1}+(n^{-1}-n^{-2}))\\
&\geq 1.&&\qedhere
\end{align*}

\end{proof}

\begin{lem}
\label{lem:other_shitty_slud_shitty}
When 
$\frac 1 2<p<\frac 1 2 + \frac 1 2 \left(1-\frac{4(\sqrt{n(n+1)} - 1)}{4n+2}\right)^{1/2}$
and $n$ is odd,
\[
\frac {2p-1}{\sqrt{n+1}} \leq \sqrt{-(n+1)\ln(4\si^2)} - 
\frac {np - (n+1)/2}{\si\sqrt{n}}.
\]
\end{lem}
\begin{proof}
The conditions on $p$ are equivalent to
\[
0 \geq (4n+2)p^2 - (4n+2)p + \sqrt{n(n+1)} - 1 = -4n\si^2 -2\si^2 + \sqrt{n(n+1)} - 1.
\]
Note that $y^{-1/2}$ is lower bounded by its Taylor expansion $1 + (1-y)/2$, meaning
\[
1+2\si = 1 + (1 - (1-4\si^2))/2 \leq (1-4\si^2)^{-1/2} = (1-2p)^{-1}.
\]
Furthermore $4\si^2 \leq 2\si$ (since $\sqrt\cdot$ is nondecreasing on $[0,1]$), so 
the above condition implies
\[
\sqrt{n(n+1)} \leq 2n\si + (1-2p)^{-1}.
\]
With some rearranging, this becomes
\[
\frac {2p-1}{\sqrt{n+1}} \leq \sqrt{(n+1)(1-4\si^2)}  - \frac {np - (n+1)/2}{\si\sqrt n},
\]
and since $-\ln(4\si^2) \geq 1-4\si^2$, the statement follows.
\end{proof}

\begin{lem}
\label{lem:slud_fuckballs}
When $m$ is even and $p\in[\bar p,1/2]$,
\[
\frac {e^{l(m)}(3/2-p)(2\si)^m}{\sqrt{2\pi m}} \geq (1 - 2e^{l(m)}\si)\min\{
\frac {\phi(\beta)}{\beta}
,\frac 1 2
\},
\]
where $\beta = ((1/2 - p)\sqrt{m})/\si$
and $\bar p= \sqrt{3/4 - 1/(2e^{l(m)})}$.
\end{lem}
\begin{proof}
First, notice that the left hand side is decreasing along $p\in[p,1/2$. Starting
from the derivative:
\begin{align*}
&\frac {d}{dp}\left(
\left(\frac {e^{l(m)}}{\sqrt{2\pi m}}\right)(\frac 3 2(p(1-p))^{m/2} - p(p(1-p))^{m/2})
\right)\\
&=
\left(\frac {e^{l(m)}}{\sqrt{2\pi m}}\right)\Bigg(\frac {3m}{4}(p(1-p))^{m/2-1}(1-2p) \\
&\quad\quad\quad\quad- (p(1-p))^{m/2}
- \frac {mp}{2}(p(1-p))^{m/2-1}(1-2p)
\Bigg)\\
&=
\left(\frac {e^{l(m)}}{\sqrt{2\pi m}}\right)
(p^2(1+m) - p(1+2m) + \frac {3m} 4).
\end{align*}
Setting this to zero and solving the quadratic yields
\[
p = \frac {1+2m \pm \sqrt{1+m + m^2 }}{2 + 2m}.
\]
Notice that since $(1+2m)/(2+2m) \geq (1+m)/(2+2m) = 1/2$, at least one of the
solutions exceeds $1/2$.  Let $p^*$ denote the solution subtracting
the discriminant.  Below, it is shown that $p^*\leq \bar p$.  But since the derivative 
evaluated at $1/2$ is negative, it must follow that the left hand side is decreasing
along $[\bar p, 1/2]$.

Since $(1/2+m)^2 = 1/4 + m + m^2 \leq 1+m+m^2$, it follows that
\[
p^* \leq \frac {1+2m - (1/2+m)}{2+2m} = \frac {1+2m}{4+4m}.
\]
By inspection, $p^* < \bar p$ for $m\in \{2,4,\ldots,28\}$.  
Next, since $-l(m) \leq (3m)^{-1}$,
\[
\bar p = \sqrt{3/4 - e^{-l(m)}/2}
\geq \sqrt{3/4 - e^{-1/3m}/2}.
\]
Any convex function may be upper bounded by its secant along an interval; thus
along $[0,1/6]$, $e^x \leq 1 + 6x(e^{1/6}-1)$.  Furthermore, $(3m)^{-1} \in [0,1/6]$ since
$m\geq 2$, thus
\begin{equation}
\label{fuckballs_dickballs_argh}
\bar p \geq \sqrt{\frac 3 4 - \frac 1 {2} \left(1 + \frac {e^{1/6}-1}{1/6} \frac 1 {3m}\right) }
= \sqrt{1/4 - \frac {e^{1/6}-1}{m} }.
\end{equation}
Next, note that
\begin{align*}
&m\bar p^2 (4+4m)^2 - m(1+2m)^2 \\
&= 4m(1+m)^2 - 16(e^{1/6}-1)(1+2m)^2 - m(1+2m)^2 \\
&= 3m + 4m^2 - 16(e^{1/6}-1)(1+2m)^2 - m(1+2m)^2,
\end{align*}
which exceeds zero when $m\geq 30$.  Rearranging, this yields $\bar p \geq p^*$.

Now consider the right hand side of the lemma statement, note that it is upper bounded by
$1/2 - e^{l(m)}\si$, which decreases as $p\to 1/2$ since $\si$ increases along this interval.
Combining all these pieces, to prove the inequality, it suffices to show that the left hand
side at $1/2$ exceeds the right hand side at $\bar p$.

First upper bound the quantity on the right hand side.  Continuing from \ref{fuckballs_dickballs_argh},
and using the secant approximation $\sqrt{x} \geq 1/6 + 4x/3$ along $[1/16,1/4]$,
\[
\bar p
\geq \sqrt{1/4 - \frac {e^{1/6}-1}{m} }
\geq \frac 1 6 + \frac 4 3\left(1/4 - \frac {e^{1/6}-1}{m}\right)
\geq \frac 1 2 - \frac 1 {4m}.
\]
On the other hand, since $\sqrt{x} \leq x + 1/4$ (first-order Taylor (tangent) at $1/4$),
$e^x\geq 1+x$, 
and $(4m)^{-1}\leq -l(m)$,
\begin{align*}
\sqrt{3/4 - e^{-l(m)}/2} 
&\leq \sqrt{3/4 - e^{-1/4m}/2}
\leq \sqrt{3/4 - (1-1/4m)/2}\\
&= \sqrt{1/4 - 1/8m}
\leq 1/2 - 1/8m.
\end{align*}
Combining these,
\[
2\si = \sqrt{4p(1-p)} \geq \sqrt{1 - 1/4m - 1/8m^2}\geq 1 - 1/4m - 1/8m^2\geq 1-1/3m
\]
the last two steps following since $\sqrt{x} \geq x$ on $[0,1]$, and $-1/8m^2 \geq -1/12m$.
Inserting this into the right hand side (with $1/2$ in the $\min$) and using
$e^{l(m)} \geq 1 - 1/3m$,
\[
\frac 1 2 (1 -2e^{l(m)}\si)
\leq \frac 1 2 (1 -(1-1/3m)^2)
\leq 1/3m 
.
\]
On the other hand, the left hand side (at $p=1/2$) may be lower bounded as
\[
\frac {e^{l(m)}}{\sqrt{2\pi m}}
\geq \frac {1-1/3m}{\sqrt{2\pi m}}
\geq \frac {5/6}{\sqrt{2\pi m}}.
\]
Comparing the square of this lower bound on the left hand side, and the square on
the upper bound of the right hand side, it is clear the left hand side is
greater, thus completing the proof.
\end{proof}

\section{Comparison of Existing Central Tail Upper Bounds}
\label{sec:pf_chern_cmp}
Applying Hoeffding, Bernstein, and Chernoff bounds to 
$\P[B(p,n)\geq n/2]$ with even $m$ and $p<1/2$ yields
\begin{align}
\P[B(p,m)\geq m/2]&\leq \exp[-m(1-2p)^2/2] &\textrm{Hoeffding,} \notag\\
\P[B(p,m)\geq m/2]&\leq \exp\left[-\frac{3m(1/2-p)^2}{(1+4p)(1-p)}\right] &\textrm{Bernstein,} \notag\\
\P[B(p,m)\geq m/2]&\leq [4p(1-p)]^{m/2} = \exp(\frac m 2 \ln(4p(1-p)))&\textrm{Chernoff.} \notag
\end{align}
\begin{thm}
The following inequalities relate the performance of these bounds:
\begin{align*}
[4p(1-p)]^{m/2} 
&\leq \exp\left[-\frac{3m(1/2-p)^2}{(1+4p)(1-p)}\right] 
\leq \exp[-m(1-2p)^2/2] &\textrm{if $p\in[0,1/4]$}; \\ 
[4p(1-p)]^{m/2} 
&\leq \exp[-m(1-2p)^2/2] 
\leq \exp\left[-\frac{3m(1/2-p)^2}{(1+4p)(1-p)}\right] 
&\textrm{if $p\in[1/4,1/2]$}. 
\end{align*}
\end{thm}
\begin{proof}
Since $\exp$ is monotonic increasing, bounds will be compared by looking solely at the
exponents.

First, note that the Chernoff Bound is always better than the Hoeffding bound:
\begin{align*}
\frac m 2 \ln(4p(1-p))
&\leq \frac m 2 (4p(1-p)-1)
= -\frac m 2 (1-2p)^2.
\end{align*}
Next, the ratio of Hoeffding's bound to Bernstein's bound is
\[
\frac 2 3 (1+4p)(1-p).
\]
By setting this quadratic to 1 and solving, the ratio is equal when $p\in\{1/4, 1/2\}$.
Furthermore, it is concave, and attains a maximum at $p = 3/8$.  Combining these facts,
it must be the case that the ratio is less than or equal to 1 along $p\in[0,1/4]$, and greater than or
equal to 1 along $p\in[1/4,1/2]$.

To finish, it remains to be shown that the Chernoff bound is less than the Bernstein bound
along $p\in[0,1/4]$.  Let $f,g$ denote the Chernoff and Bernstein bounds, respectively.
$df/dp\geq 0$ and $dg/dp\geq 0$  along this interval, meaning both are increasing.  Thus the
result follows by the fact that $f(1/8) < g(0)$ and $f(1/4) < g(1/8)$.
\end{proof}

\section{Riemann Sums of Convex, Decreasing Functions}
\label{sec:riemann_lemma}
\begin{figure}
\includegraphics[width=0.50\textwidth]{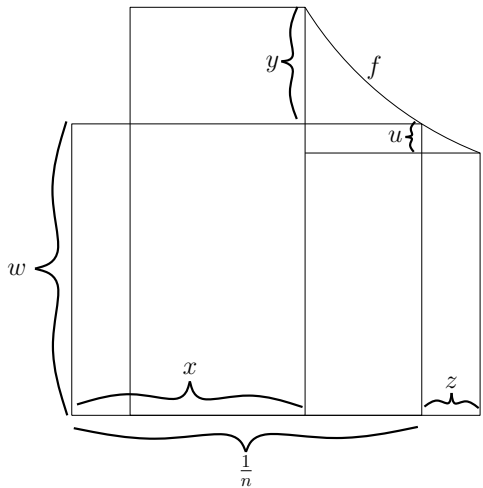}
\caption{Given a convex, decreasing function $f$, forming Riemann sums
with successively finer subdivisions yields a sequence of monotonic increasing
approximants.}
\label{fig:psi_riemann}
\end{figure}
\begin{lem}
Let intervals $(a,b]\subseteq (c,d]$, integers $n < m$ with $(d-c)/m \leq (b-a)/n$, and a function $f$ 
convex and decreasing on $(c,d]$ be given.  Let $R(f;(t,u],v)$ denote the Riemann
sum approximation of $\int_t^u f$ consisting of $v$ equal pieces, whose height
is $f$ evaluated at their right endpoint.  Then $R(f;(a,b],n) \leq R(f;(c,d],m)$.
\end{lem}
\begin{proof}
Let $I$ denote some block of $R(f;(a,b],n)$; that is, there is some integer
$1\leq j \leq n$ such that $I = [a + j(b-a)/m, a + (j+1)(b-a)/m]$.  The desired 
fact will be shown for the restriction to $I$, from which the general statement follows
by considering all $n$ pieces of $[a,b]$.

The ingredients of the proof appear in figure~\ref{fig:psi_riemann}.
In particular, consider the greatest point within $I$ where the right endpoint
of a block of $R(f;(c,d],m)$ falls (such a point must exist since $R(f;(c,d],m)$ has
narrower blocks).  Let $w= f(a+(j+1)(b-a)/m)$ be the height of the block, and let 
$x$ be the offset at which this point falls.  Furthermore, let $y$ be the height
above $w$ of the block.

Now consider the immediately next block (in $R(f;(c,d],m)$); its boundary must
fall outside $I$.  Let $z$ denote this distance beyond the edge of $I$, and let
$u$ be the height below $w$ of this point.  By assumption, it must be the case
that $z\leq x$.  Furthermore, since $f$ is convex on $(c,d]$, it must be the case 
that $\frac {y}{1/n-x} \geq \frac u z$.  Thus, to compute the area with respect 
to $R(f;(c,d],m)$,
one has the lower bound
\[
x(w+y) + \left(\frac 1 n - x\right)(w-u)
\geq  \frac w n +xy - \left(\frac 1 n - x\right) \frac {yz}{1/n-x}
\geq \frac w n.
\]
Note that this last quantity is precisely the area of the block with respect
to $R(f;(a,b],n)$, which
completes the proof.
\end{proof}

\end{document}